\documentclass[11pt]{article}

\usepackage{amssymb,latexsym}
\usepackage{epsfig}
\usepackage{eufrak}
\usepackage{amsmath}
\usepackage{mathrsfs}
\usepackage{color}

\newtheorem{theorem}{Theorem}
\newtheorem{lemma}{Lemma}

\newcommand{\be}{\begin{equation}}
\newcommand{\ee}{\end{equation}}
\newcommand{\bee}{\begin{eqnarray*}}
\newcommand{\eee}{\end{eqnarray*}}
\newcommand{\bel}{\begin{eqnarray}}
\newcommand{\eel}{\end{eqnarray}}
\newcommand{\bec}{\begin{cases}}
\newcommand{\eec}{\end{cases}}
\newcommand{\bem}{\begin{bmatrix}}
\newcommand{\eem}{\end{bmatrix}}

\newcommand{\la}{\label}
\newcommand{\li}{\left}
\newcommand{\ri}{\right}

\newcommand{\ovl}{\overline}

\newcommand{\vep}{\varepsilon}
\newcommand{\lm}{\lambda}

\newcommand{\si}{\sigma}

\newcommand{\de}{\delta}

\newcommand{\vDe}{\varDelta}

\newcommand{\vse}{\vartheta}
\newcommand{\se}{\theta}

\newcommand{\Om}{\Omega}

\newcommand{\f}{\frac}

\newcommand{\cd}{\cdots}

\newcommand{\qu}{\quad}
\newcommand{\qqu}{\qquad}
\newcommand{\fa}{\forall}

\newcommand{\mscr}{\mathscr}
\newcommand{\mcal}{\mathcal}

\newcommand{\bb}{\mathbb}

\newcommand{\wh}{\widehat}

\newcommand{\mrm}{\mathrm}
\newcommand{\bs}{\boldsymbol}

\newcommand{\tx}{\text}

\newcommand{\iy}{\infty}

\newcommand{\pa}{\partial}

\newcommand{\bed}{\begin{description}}
\newcommand{\eed}{\end{description}}
\newcommand{\bei}{\begin{itemize}}
\newcommand{\eei}{\end{itemize}}
\newcommand{\ben}{\begin{enumerate}}
\newcommand{\een}{\end{enumerate}}
\newcommand{\bib}{\bibitem}
\newcommand{\beL}{\begin{lemma}}
\newcommand{\eeL}{\end{lemma}}
\newcommand{\beT}{\begin{theorem}}
\newcommand{\eeT}{\end{theorem}}
\newcommand{\sect}{\section}

\newcommand{\bpf}{\begin{pf}}
\newcommand{\epf}{\end{pf}}
\newcommand{\bsk}{\bigskip}

\newcommand{\pfbox}{\hfill\mbox{$\Box$}}

\newenvironment{pf}{\paragraph*{Proof{\rm.}}}{\pfbox\bigskip}

\begin{document}

\title{{\bf On Estimation and Optimization of Mean Values of Bounded Variables}
\thanks{The author is currently with Department of Electrical Engineering,
Louisiana State University at Baton Rouge, LA 70803, USA, and Department of Electrical Engineering, Southern University and A\&M College, Baton
Rouge, LA 70813, USA; Email: chenxinjia@gmail.com}}

\author{Xinjia Chen}

\date{First Submitted in February  2008}

\maketitle

\begin{abstract}

In this paper, we develop a general approach for probabilistic estimation and optimization. An explicit formula and a computational approach are
established for controlling the reliability of probabilistic estimation based on a mixed criterion of absolute and relative errors. By employing
the Chernoff-Hoeffding bound and the concept of sampling, the minimization of a probabilistic function is transformed into an optimization
problem amenable for gradient descendent algorithms.

\end{abstract}

\section{Analytical Sample Size Formula for Estimation of Mean Values}

Let $X$ be a random variable bounded in interval $[0, 1]$ with mean $\bb{E} [X] = \mu \in (0, 1)$, which are defined on a probability space
$(\Om, \mscr{F}, \Pr)$.  In many areas of sciences and engineering, it is desired to estimate $\mu$ based on samples $X_1, X_2, \cd, X_n$ of
$X$. Frequently,  the samples $X_1, X_2, \cd, X_n$ may not be identical and independent (i.i.d). Thus, it is a significant problem to estimate
$\mu$ under the assumption that \bel &  & 0 \leq X_k \leq 1 \qu \tx{almost surely for any positive integer $k$},  \la{gen89b}\\
&  & \bb{E} [ X_k \mid \mscr{F}_{k-1} ] = \mu  \qu \tx{almost surely for any positive integer $k$},  \la{gen89} \eel where  $\{ \mscr{F}_k, \; k
= 0, 1, \cd, \iy \}$ is a sequence of $\si$-subalgebra such that $\{ \emptyset, \Om \} = \mscr{F}_0 \subset \mscr{F}_1 \subset \mscr{F}_2
\subset \cd \subset \mscr{F}$, with $\mscr{F}_k$ being generated by $X_1, \cd, X_k$.

Naturally, an estimator for $\mu$ is taken as \be \la{est}
 \wh{\bs{\mu}} = \f{\sum_{i = 1}^n X_i } { n
 }.
 \ee
    Since $\wh{\bs{\mu}}$ is of random nature, it is crucial to control the statistical error. For this purpose, we have established the
 following result.

\beT Let $\de \in (0,1)$.  Let $\vep_a \in (0, 1)$ and $\vep_r \in (0, 1)$ be real numbers such that $\f{ \vep_a}{ \vep_r} + \vep_a \leq
\f{1}{2}$.  Assume that (\ref{gen89b}) and (\ref{gen89}) are true.  Then,
 \be
 \la{cov}
 \Pr \li \{ |\wh{\bs{\mu}} - \mu| < \vep_a  \; \; \mrm{or} \;\; \li | \f{
\wh{\bs{\mu}} - \mu } {\mu} \ri | < \vep_r  \ri \} > 1 - \de \ee
 for any $\mu \in (0, 1)$ provided that \be \la{con}
 n > \f{ \vep_r \ln \f{2}{\de} } { \li ( \vep_a + \vep_a \vep_r \ri ) \ln (1 + \vep_r) + \li ( \vep_r - \vep_a -
\vep_a \vep_r \ri ) \ln \li ( 1 -  \f{ \vep_a \vep_r } { \vep_r - \vep_a  }
 \ri )  }.
\ee \eeT

\bsk

It should be noted that conventional methods for determining sample sizes are based on normal approximation, see \cite{Desu} and the references
therein.  In contrast, Theorem 1 offers a rigorous method for determining sample sizes.  In the special case that $X$ is a Bernoulli random
variable, a numerical approach has been developed by Chen \cite{Chen} which permits exact computation of the minimum sample size.

\section{A Computational Approach for General Case}

In this section, we shall investigate an exact computational sample size method for the case that $X \in [a, b]$ with $\bb{E} [ X] = \mu$.
Assume that \bel &  & a \leq X_k \leq b \qu \tx{almost surely for any positive integer $k$},  \la{gen89b88a}\\
&  & \bb{E} [ X_k \mid \mscr{F}_{k-1} ] = \mu  \qu \tx{almost surely for any positive integer $k$},  \la{gen8988b} \eel where  $\{ \mscr{F}_k,
\; k = 0, 1, \cd, \iy \}$ is a sequence of $\si$-subalgebra such that $\{ \emptyset, \Om \} = \mscr{F}_0 \subset \mscr{F}_1 \subset \mscr{F}_2
\subset \cd \subset \mscr{F}$, with $\mscr{F}_k$ being generated by $X_1, \cd, X_k$.

We wish to determine minimum sample size $n$ such that \be \la{wish}
 \Pr \li \{ |\wh{\bs{\mu}} - \mu| < \vep_a  \; \; \mrm{or} \;\; \li |
\wh{\bs{\mu}} - \mu
 \ri | < \vep_r  | \mu | \ri \} > 1 - \de
\ee for any $\mu \in [a, b]$, where $\wh{\bs{\mu}}$ is defined by (\ref{est}).  Unlike the special case that $X$ is bounded in interval $[0,
1]$, there is no explicit formula for the general case that $X$ is bounded in interval $[a, b]$.  We will employ the branch and bound technique
of global optimization. For this purpose, we need to derive a sample size formula and the associated bounding method.

To describe the relevant theory for computing sample sizes, define function
\[
\mscr{M} (z,\se) = \bec z \ln \f{\se}{z} + (1 - z) \ln \f{1 - \se}{1 - z} &
\tx{for} \; z \in (0,1) \; \tx{and} \; \se \in (0, 1),\\
\ln(1-\se) & \tx{for} \; z = 0 \; \tx{and} \; \se \in (0, 1),\\
\ln \se &  \tx{for} \; z = 1 \; \tx{and} \; \se \in (0, 1),\\
- \iy &  \tx{for} \; z \in [0, 1] \; \tx{and} \; \se \notin (0, 1) \eec
\]
Define {\small \bee &  & \vse(\mu) = \f{\mu -
a} {b - a},\\
& &  g(\mu) = \vse(\mu) -  \f{\max \{ \vep_a, \; \vep_r | \mu | \} }{b - a}, \\
& & h(\mu) = \vse(\mu) +  \f{\max \{ \vep_a, \; \vep_r | \mu | \} }{b - a},\\
&  & \mcal{W}(\mu) = \max \li \{ \mscr{M} \li ( g(\mu), \vse(\mu) \ri ), \; \mscr{M} \li ( h(\mu), \vse(\mu) \ri ) \ri \} \eee} for $\mu \in [a,
b]$. By virtue of such functions, we have established theoretical results which are essential for the exact computation of sample sizes as
follows.

\beT \la{Bounded_Mean_mix_general_Massart} Assume that  (\ref{gen89b88a}) and (\ref{gen8988b}) are satisfied. Then, (\ref{wish}) holds for any
$\mu \in [a, b]$ provided that {\small \be \la{need98} n \geq \f{\ln \f{\de}{2} } { \max_{\nu \in [a, b]} \mcal{W} (\nu) }.  \ee} Moreover, \bel
& & \mcal{W} (\nu) \leq \max \li \{ \mscr{M} \li ( g(d) , \vse (c) \ri ),
\; \mscr{M} \li ( h(c), \vse (d) \ri ) \ri \},  \la{bounda}\\
&  & \mcal{W} (\nu) \geq \max \li \{ \mscr{M} \li ( g(c), \vse (d) \ri ), \; \mscr{M} \li ( h(d), \vse (c) \ri ) \ri \} \la{boundb} \eel for
$\nu \in [c, d] \subseteq [a, b]$ such that $g(d) \leq \vse (c) \leq \vse(d) \leq h(c)$. \eeT

See Appendix \ref{Bounded_Mean_mix_general_Massart_app} for a proof.

Since  (\ref{bounda}) and (\ref{boundb}) of Theorem \ref{Bounded_Mean_mix_general_Massart} provide computable upper and lower bounds of
 $\mcal{W} (\nu)$, the  maximum of $\mcal{W} (\nu)$ over $[a, b]$ can be exactly computed with the  {\it Branch and Bound} method proposed by Land and Doig
\cite{Land}.

\sect{Optimization of Probability}

In many applications, it is desirable to find a vector of real numbers $\se$ to minimize a probability, $p(\se)$, which can be expressed as
\[
p(\se) = \Pr \{ Y (\se, \bs{\vDe}) \leq 0 \},
\]
where $Y(\se, \bs{\vDe})$ is piece-wise continuous with respect to $\se$ and $\bs{\vDe}$ is a random vector.  If we define
\[
\mu (\lm, \se) = \bb{E} [ e^{ - \lm Y (\se, \bs{\vDe})  } ],
\]
then, applying Chernoff bound \cite{Chernoff}, we have
 \[
p(\se) \leq \inf_{\lm > 0} \mu (\lm, \se).
 \]
This indicates that we can make $p(\se)$ small by making $\mu (\lm, \se)$ small.  Hence, we shall attempt to minimize $\mu (\lm, \se)$ with
respect to $\lm
> 0$ and $\se$.

To make the new objective function $\mu (\lm, \se)$ more tractable, we take a sampling approach. Specifically, we obtain $n$ i.i.d. samples
$\bs{\vDe}_1, \cd, \bs{\vDe}_n$ of $\bs{\vDe}$ and approximate $\mu (\lm, \se) $ as
\[
g(\lm, \se) = \f{ \sum_{i=1}^n e^{ - \lm Y (\se, \bs{\vDe_i})  }  } { n }.
\]
A critical step is the determination of sample size $n$ so that $g(\lm, \se)$ is sufficiently close to $\mu(\lm, \se)$. Since $0 < e^{ - \lm Y
(\se, \bs{\vDe}) }
 < 1$, an appropriate value of $n$ can be computed based on (\ref{con}) of Theorem 1.

Finally, we have transformed the problem of minimizing the probability function $p(\se)$ as the problem of minimizing a piece-wise continuous
function $g(\lm, \se)$.  Since $g(\lm, \se)$ is a more smooth function, we can bring all the power of nonlinear programming to solve the
problem. An extremely useful tool is the {\it gradient descendent algorithm}, see, e.g. \cite{ba} and the references therein.

\sect{Proof of Theorem 1}

To prove the theorem, we shall introduce function
\[
\psi (\vep, \mu) =  (\mu + \vep) \ln \f{ \mu}{\mu + \vep} + (1 - \mu - \vep) \ln \f{1 - \mu}{1 - \mu - \vep}
\]
where $0 < \vep < 1 - \mu$.   We need some preliminary results.

The following lemma is due to Hoeffding \cite{Hoeffding}.

\beL \la{lem1} Assume that (\ref{gen89b}) and (\ref{gen89}) hold for any positive integer $k$. Then,
\[
\Pr \{ \wh{\bs{\mu}} \geq \mu + \vep \} \leq \exp ( n \; \psi (\vep, \mu) ) \qu \tx{for} \qu 0 < \vep < 1 - \mu < 1,
\]
\[
\Pr \{ \wh{\bs{\mu}} \leq \mu - \vep \} \leq \exp ( n \; \psi (-\vep, \mu) ) \qu \tx{for} \qu 0 < \vep < \mu < 1.
\]
\eeL

\beL  \la{lem2}

Let $0 < \vep < \f{1}{2}$. Then, $\psi (\vep, \mu)$ is monotonically increasing with respective to $\mu \in (0,  \f{1}{2} - \vep)$ and
monotonically decreasing with respective to $\mu \in (\f{1}{2},  1 - \vep)$.  Similarly, $\psi (- \vep, \mu)$ is monotonically increasing with
respective to $\mu \in (\vep, \f{1}{2})$ and monotonically decreasing with respective to $\mu \in (\f{1}{2} + \vep,  1)$.\eeL

\bpf

Tedious computation shows that
\[
\f{ \pa  \psi (\vep, \mu) } { \pa \mu } = \ln \f{ \mu ( 1 - \mu - \vep) } { (\mu + \vep) (1 - \mu) }  + \f{ \vep } { \mu} + \f{ \vep } { 1 -
\mu}
\]
and
\[
\f{ \pa^2  \psi (\vep, \mu) } { \pa \mu^2 } = - \f{ \vep^2 } { \mu^2 (\mu+\vep) } - \f{ \vep^2 } { (1 - \mu)^2 (1 - \mu - \vep) } < 0
\]
for $0 < \vep < 1 - \mu < 1$.  Note that
\[
\f{ \pa  \psi (\vep, \mu) } { \pa \mu } |_{\mu = \f{1}{2}} =  \ln \f{ 1 - 2 \vep  } { 1 + 2 \vep } + \vep < 0
\]
because
\[
\f{d \li [  \ln \f{ 1 - 2 \vep  } { 1 + 2 \vep } + \vep \ri ]} { d \vep } = - \f{4}{1 - 4 \vep^2} < 0.
\]
Moreover,
\[
\f{ \pa  \psi (\vep, \mu) } { \pa \mu } |_{\mu = \f{1}{2} - \vep }  = \ln \f{ 1 - 2 \vep  } { 1 + 2 \vep } + \f{4 \vep}{1 - 4 \vep^2} > 0
\]
because
\[
\f{ d \li [ \ln \f{ 1 - 2 \vep  } { 1 + 2 \vep } + \f{4 \vep}{1 - 4 \vep^2} \ri ] }{d \vep  } = \f{ 32 \vep^2 } { (1 - \vep^2)^2} > 0.
\]

Similarly,
\[
\f{ \pa  \psi (- \vep, \mu) } { \pa \mu } = \ln \f{ \mu ( 1 - \mu + \vep) } { (\mu - \vep) (1 - \mu) }  - \f{ \vep } { \mu} - \f{ \vep } { 1 -
\mu}
\]
and
\[
\f{ \pa^2  \psi (- \vep, \mu) } { \pa \mu^2 } = - \f{ \vep^2 } { \mu^2 (\mu - \vep) } - \f{ \vep^2 } { (1 - \mu)^2 (1 - \mu + \vep) } < 0
\]
for $0 < \vep < \mu < 1$.  Hence,
\[
\f{ \pa  \psi (- \vep, \mu) } { \pa \mu } |_{\mu = \f{1}{2}} =  \ln \f{ 1 + 2 \vep  } { 1 - 2 \vep } - \vep > 0
\]
because
\[
\f{d \li [  \ln \f{ 1 + 2 \vep  } { 1 - 2 \vep } - \vep \ri ]} { d \vep } =  \f{4}{1 - 4 \vep^2} > 0;
\]
and
\[
\f{ \pa  \psi (- \vep, \mu) } { \pa \mu } |_{\mu = \f{1}{2} + \vep }  = \ln \f{ 1 + 2 \vep  } { 1 - 2 \vep } - \f{4 \vep}{1 - 4 \vep^2} < 0
\]
as a result of
\[
\f{ d \li [ \ln \f{ 1 + 2 \vep  } { 1 - 2 \vep } - \f{4 \vep}{1 - 4 \vep^2} \ri ] }{d \vep  } = - \f{ 32 \vep^2 } { (1 - \vep^2)^2} < 0.
\]

Since $\f{ \pa  \psi (\vep, \mu) } { \pa \mu } |_{\mu = \f{1}{2}} < 0, \; \f{ \pa  \psi (\vep, \mu) } { \pa \mu } |_{\mu = \f{1}{2} - \vep }
> 0$ and $\psi (\vep, \mu)$ is concave with respect to $\mu$, it must be true that $\psi (\vep, \mu)$ is monotonically increasing with respective to $\mu \in (0,  \f{1}{2} - \vep)$ and
monotonically decreasing with respective to $\mu \in (\f{1}{2},  1 - \vep)$. Since $\f{ \pa  \psi (- \vep, \mu) } { \pa \mu } |_{\mu = \f{1}{2}}
> 0, \; \f{ \pa  \psi (- \vep, \mu) } { \pa \mu } |_{\mu = \f{1}{2} + \vep } < 0$ and $\psi (\vep, \mu)$ is concave with respect to $\mu$, it must be
true that $\psi (- \vep, \mu)$ is monotonically increasing with respective to $\mu \in (\vep, \f{1}{2})$ and monotonically decreasing with
respective to $\mu \in (\f{1}{2} + \vep,  1)$.

\epf

\beL \la{lem3}

Let $0 < \vep < \f{1}{2}$.  Then,
\[
\psi (\vep, \mu) > \psi (-\vep, \mu) \qqu \fa \mu \in \li (\vep, \f{1}{2} \ri ],
\]
\[
\psi (\vep, \mu) < \psi (-\vep, \mu) \qqu \fa \mu \in \li (\f{1}{2}, 1 - \vep \ri ).
\]
\eeL

\bpf

It can be shown that
\[
\f{ \pa  [\psi (\vep, \mu) - \psi (-\vep, \mu)] } { \pa \vep } = \ln \li [ 1 + \f{ \vep^2 ( 1 - 2 \mu) } { (\mu^2 - \vep^2) ( 1 - \mu)^2 } \ri ]
\]
for $0 < \vep < \min (\mu, 1 - \mu)$.  Note that \[ \f{ \vep^2 ( 1 - 2 \mu) } { (\mu^2 - \vep^2) ( 1 - \mu)^2 }  > 0 \qu \mrm{for} \qu  \vep <
\mu < \f{1}{2}
\]
and \[ \f{ \vep^2 ( 1 - 2 \mu) } { (\mu^2 - \vep^2) ( 1 - \mu)^2 }  < 0 \qu \mrm{for} \qu  \vep < \f{1}{2} < \mu < 1 - \vep.
\]
Therefore,
\[ \f{ \pa  [\psi (\vep, \mu) - \psi (-\vep, \mu)] } { \pa \vep }  > 0 \qu \mrm{for} \qu  \vep < \mu < \f{1}{2}
\]
and \[ \f{ \pa  [\psi (\vep, \mu) - \psi (-\vep, \mu)] } { \pa \vep } < 0 \qu \mrm{for} \qu  \vep < \f{1}{2} < \mu < 1 - \vep.
\]
So, we can complete the proof of the lemma by observing the sign of the partial derivative {\small $\f{ \pa  [\psi (\vep, \mu) - \psi (-\vep,
\mu)] } { \pa \vep }$} and the fact that $\psi (\vep, \mu) - \psi (-\vep, \mu) = 0$ for $\vep = 0$.

\epf

\beL \la{lem4}

Let $0 < \vep < 1$.  Then, $\psi  \li ( \vep \mu, \mu \ri )$ is monotonically decreasing with respect to {\small $\mu \in \li (0, \f{1}{1 +
\vep} \ri )$}. Similarly, $\psi  \li (- \vep \mu, \mu \ri )$ is monotonically decreasing with respect to $\mu \in (0, 1)$. \eeL

\bpf Note that
\[
\f{ \pa   \psi  \li ( \vep \mu, \mu \ri )} { \pa \mu } = (1 + \vep) \ln \f{ 1 - (1 + \vep) \mu } { 1 - \mu } - (1 + \vep) \ln (1 + \vep) +
\f{\vep}{1 - \mu}
\]
and
\[
\f{ \pa^2   \psi  \li ( \vep \mu, \mu \ri )} { \pa \mu^2 } = - \f{ \vep^2 } { (1 - \mu)^2 [1 - (1 + \vep) \mu ] } < 0
\]
for any $\mu \in \li (0, \f{1}{1 + \vep} \ri )$.

Since $\f{ \pa   \psi  \li ( \vep \mu, \mu \ri )} { \pa \mu }|_{\mu = 0} = \vep  - (1 + \vep) \ln (1 + \vep) < 0$, we have
\[
\f{ \pa   \psi  \li ( \vep \mu, \mu \ri )} { \pa \mu } < 0, \qu \fa \mu \in \li (0, \f{1}{1 + \vep} \ri )
\]
and it follows that $\psi  \li ( \vep \mu, \mu \ri )$ is monotonically decreasing with respect to {\small $\mu \in \li (0, \f{1}{1 + \vep} \ri
)$}.

Similarly, since
\[
\f{ \pa   \psi  \li ( - \vep \mu, \mu \ri )} { \pa \mu }|_{\mu = 0} = - \vep  - (1 - \vep) \ln (1 - \vep) < 0 \] and \[ \f{ \pa^2   \psi  \li (
\vep \mu, \mu \ri )} { \pa \mu^2 } = - \f{ \vep^2 } { (1 - \mu)^2 [1 - (1 - \vep) \mu ] } < 0, \qu \fa \mu \in (0, 1)
\]
 we have
\[
\f{ \pa   \psi  \li ( - \vep \mu, \mu \ri )} { \pa \mu } < 0, \qu \fa \mu \in (0, 1)
\]
and, consequently, $\psi  \li (- \vep \mu, \mu \ri )$ is monotonically decreasing with respect to $\mu \in (0, 1)$.

\epf

\beL \la{lem5}  Suppose $0 < \vep_r < 1$ and $0 < \f{\vep_a}{\vep_r} + \vep_a \leq \f{1}{2}$.  Then, \be \la{ineqm}
 \Pr \{ \wh{\bs{\mu}} \leq \mu - \vep_a \} \leq \exp \li (n \; \psi  \li ( -
\vep_a, \f{\vep_a}{\vep_r}\ri ) \ri ) \ee for $0 < \mu \leq \f{\vep_a}{\vep_r}$.  \eeL

\bpf We shall show (\ref{ineqm}) by investigating three cases as follows.  In the case of $\mu < \vep_a$, it is clear that
\[
\Pr \{ \wh{\bs{\mu}} \leq \mu - \vep_a \} = 0  < \exp \li (n \; \psi  \li ( - \vep_a, \f{\vep_a}{\vep_r}\ri ) \ri ).
\]

In the case of $\mu = \vep_a$, we have \bee \Pr \{ \wh{\bs{\mu}} \leq \mu - \vep_a \} & = & \lim_{\eta \uparrow \vep_a} \Pr \{ \wh{\bs{\mu}}
\leq \mu - \eta \}\\
&  \leq & \lim_{\eta \uparrow \vep_a} \exp \li (n \; \psi  \li ( - \eta, \mu \ri ) \ri ) = \exp \li (n \; \psi  \li ( - \vep_a, \mu \ri ) \ri )\\
& = &  \exp \li (n \; \psi  \li ( - \vep_a, \vep_a \ri ) \ri )\\
& < &  \exp \li (n \; \psi  \li ( - \vep_a, \f{\vep_a}{\vep_r}\ri ) \ri ), \eee where the last inequality follows from Lemma \ref{lem2} and the
fact that $\vep_a < \f{\vep_a}{\vep_r} \leq \f{1}{2} - \vep_a$.

In the case of $\vep_a < \mu \leq \f{\vep_a}{\vep_r}$,  we have
\[
\Pr \{ \wh{\bs{\mu}} \leq \mu - \vep_a \} \leq \exp (n \; \psi (- \vep_a, \mu)) < \exp \li (n \; \psi  \li ( - \vep_a, \f{\vep_a}{\vep_r}\ri )
\ri ),
\]
where the first inequality follows from Lemma \ref{lem1} and the second inequality follows from Lemma \ref{lem2} and the fact that $\vep_a <
\f{\vep_a}{\vep_r} \leq \f{1}{2} - \vep_a$.  So, (\ref{ineqm}) is established. \epf

\beL \la{lem6}  Suppose $0 < \vep_r < 1$ and $0 < \f{\vep_a}{\vep_r} + \vep_a \leq \f{1}{2}$.   Then, \be \la{ineqmm}
 \Pr \{ \wh{\bs{\mu}} \geq (1 + \vep_r) \mu \} \leq \exp \li (n \; \psi  \li (
\vep_a, \f{\vep_a}{\vep_r}\ri ) \ri ) \ee for $\f{\vep_a}{\vep_r} < \mu < 1$.  \eeL

\bpf We shall show (\ref{ineqmm}) by investigating three cases as follows.  In the case of $\mu > \f{1}{1 + \vep_r}$, it is clear that
\[
\Pr \{ \wh{\bs{\mu}} \geq (1 + \vep_r) \mu \} = 0  < \exp \li (n \; \psi  \li ( \vep_a, \f{\vep_a}{\vep_r}\ri ) \ri ).
\]

In the case of $\mu = \f{1}{1 + \vep_r}$, we have \bee \Pr \{ \wh{\bs{\mu}} \geq (1 + \vep_r) \mu \} & = & \lim_{\eta \uparrow \vep_r} \Pr \{
\wh{\bs{\mu}} \geq (1 + \eta) \mu \}\\
& \leq & \lim_{\eta \uparrow \vep_r} \exp (n \; \psi ( \eta \mu, \mu)) = \exp (n \; \psi ( \vep_r \mu, \mu))\\
& < &  \exp \li (n \; \psi  \li ( \vep_a, \f{\vep_a}{\vep_r}\ri ) \ri ), \eee where the last inequality follows from Lemma \ref{lem4} and the
fact that $ \f{\vep_a}{\vep_r} \leq \f{1}{2} \f{1}{1 + \vep_r} < \f{1}{1 + \vep_r}$ as a result of $0 < \f{\vep_a}{\vep_r} + \vep_a \leq
\f{1}{2}$.

In the case of $\f{\vep_a}{\vep_r} < \mu < \f{1}{1 + \vep_r}$,  we have
\[
\Pr \{ \wh{\bs{\mu}} \leq (1 + \vep_r) \mu \} \leq \exp (n \; \psi (\vep_r \mu, \mu)) < \exp \li (n \; \psi  \li ( \vep_a, \f{\vep_a}{\vep_r}\ri
) \ri ),
\]
where the first inequality follows from Lemma \ref{lem1} and the second inequality follows from Lemma \ref{lem4}.  So, (\ref{ineqmm}) is
established. \epf

\bsk

We are now in a position to prove the theorem. We shall assume (\ref{con}) is satisfied and show that (\ref{cov}) is true.  It suffices to show
that
\[
\Pr \{ |\wh{\bs{\mu}} - \mu | \geq \vep_a, \; |\wh{\bs{\mu}} - \mu | \geq \vep_r \mu \} < \de.
\]

For $0 < \mu \leq \f{\vep_a}{\vep_r}$, we have \bel \Pr \{ |\wh{\bs{\mu}} - \mu | \geq \vep_a, \; |\wh{\bs{\mu}} - \mu | \geq \vep_r \mu \} & =
& \Pr \{ |\wh{\bs{\mu}} - \mu | \geq \vep_a  \} \nonumber\\
& = & \Pr \{ \wh{\bs{\mu}}  \geq \mu + \vep_a \} + \Pr \{ \wh{\bs{\mu}} \leq \mu - \vep_a \} \la{ineq6}. \eel

Noting that $0 < \mu + \vep_a \leq \f{\vep_a}{\vep_r} + \vep_a \leq \f{1}{2}$, we have
\[
\Pr \{ \wh{\bs{\mu}}  \geq \mu + \vep_a \} \leq \exp (n \; \psi (\vep_a, \mu)) \leq \exp \li (n \; \psi  \li ( \vep_a, \f{\vep_a}{\vep_r} \ri )
\ri ),
\]
where the first inequality follows from Lemma \ref{lem1} and the second inequality follows from Lemma \ref{lem2}.  It can be checked that
(\ref{con}) is equivalent to
\[
\exp \li (n \; \psi  \li ( \vep_a, \f{\vep_a}{\vep_r} \ri ) \ri ) < \f{\de}{2}.
\]
Therefore,
\[
\Pr \{ \wh{\bs{\mu}}  \geq \mu + \vep_a \}  < \f{\de}{2}
\]
for $0 < \mu \leq \f{\vep_a}{\vep_r}$. \bsk

On the other hand, since $\vep_a < \f{\vep_a}{\vep_r} < \f{1}{2}$, by Lemma \ref{lem5} and Lemma \ref{lem3}, we have
\[
 \Pr \{ \wh{\bs{\mu}} \leq \mu - \vep_a \} \leq \exp \li (n \; \psi  \li ( -
\vep_a, \f{\vep_a}{\vep_r}\ri ) \ri ) \leq \exp \li (n \; \psi  \li ( \vep_a, \f{\vep_a}{\vep_r}\ri ) \ri ) < \f{\de}{2}
\]
for $0 < \mu \leq \f{\vep_a}{\vep_r}$.  Hence, by (\ref{ineq6}),
\[
\Pr \{ |\wh{\bs{\mu}} - \mu | \geq \vep_a, \; |\wh{\bs{\mu}} - \mu | \geq \vep_r \mu \} < \f{\de}{2} + \f{\de}{2} = \de.
\]
This proves (\ref{cov}) for $0 < \mu \leq \f{\vep_a}{\vep_r}$.

\bsk

For $\f{\vep_a}{\vep_r} < \mu < 1$, we have \bee \Pr \{ |\wh{\bs{\mu}} - \mu | \geq \vep_a, \; |\wh{\bs{\mu}} - \mu | \geq \vep_r \mu \} & = &
\Pr \{ |\wh{\bs{\mu}} - \mu | \geq \vep_r \mu \}\\
& = & \Pr \{ \wh{\bs{\mu}}  \geq \mu + \vep_r \mu \} + \Pr \{ \wh{\bs{\mu}} \leq \mu - \vep_r \mu \}. \eee  Invoking Lemma \ref{lem6}, we have
\[
\Pr \{ \wh{\bs{\mu}}  \geq \mu + \vep_r \mu \} \leq \exp \li (n \; \psi  \li ( \vep_a, \f{\vep_a}{\vep_r} \ri ) \ri ).
\]
On the other hand,
\[
 \Pr \{ \wh{\bs{\mu}}
\leq \mu - \vep_r \mu \} \leq \exp (n \; \psi ( - \vep_r \mu, \mu)) \leq \exp \li (n \; \psi  \li ( - \vep_a, \f{\vep_a}{\vep_r}\ri ) \ri ) \leq
\exp \li (n \; \psi  \li ( \vep_a, \f{\vep_a}{\vep_r}\ri ) \ri )
\]
where the first inequality follows from Lemma \ref{lem1}, the second inequality follows from Lemma \ref{lem4}, and the last inequality follows
from Lemma \ref{lem3}.  Hence,
\[
\Pr \{ |\wh{\bs{\mu}} - \mu | \geq \vep_a, \; |\wh{\bs{\mu}} - \mu | \geq \vep_r \mu \}  \leq  2 \exp \li (n \; \psi  \li ( \vep_a,
\f{\vep_a}{\vep_r} \ri ) \ri ) < \de. \] This proves (\ref{cov}) for $\f{\vep_a}{\vep_r} < \mu < 1$. The proof of Theorem 1 is thus completed.

\section{Proof of Theorem \ref{Bounded_Mean_mix_general_Massart}} \la{Bounded_Mean_mix_general_Massart_app}

Define $\ovl{Y}_n = \f{1}{n} \sum_{i = 1} Y_i$ with $Y_i = \f{ X_i - a}{b - a}$ for $i = 1, \cd, n$.  Then, $\bb{E} [ Y_i ] = \vse(\mu)$ for $i
= 1, \cd, n$.  Moreover,  \bel \Pr \{ | \ovl{X}_n - \mu | \geq \vep_a, \; | \ovl{X}_n - \mu | \geq \vep_r |\mu| \}
& = & \Pr \{   \ovl{X}_n \leq \mu - \max ( \vep_a, \vep_r |\mu | )  \} \nonumber\\
&  & \;  + \Pr \{   \ovl{X}_n \geq \mu + \max ( \vep_a, \vep_r |\mu | ) \} \nonumber\\
& = & \Pr \li \{  \ovl{Y}_n \leq g(\mu)  \ri \} + \Pr \li \{  \ovl{Y}_n \geq h(\mu)  \ri \}.   \la{firt8} \eel   It follows from (\ref{firt8})
and Lemma 1 that \bee \Pr \{ | \ovl{X}_n - \mu | \geq \vep_a, \;  | \ovl{X}_n - \mu | \geq \vep_r |\mu| \} & \leq & \exp \li (  n \mscr{M} ( g
(\mu), \vse(\mu) \ri ) + \exp \li ( n \mscr{M}
( h (\mu), \vse(\mu) \ri )\\
&  \leq & 2 \exp (n \mcal{W} (\mu) ),   \eee from which it follows immediately that (\ref{wish}) holds for any $\mu \in [a, b]$ provided that
(\ref{need98}) is true.

Now we shall show (\ref{bounda}) and (\ref{boundb}).  For $\nu \in [c, d] \subseteq [a, b]$ with $g(d) \leq \vse (c) \leq \vse(d) \leq h(c)$,
it can be shown that
\[
g(c) \leq g(\nu) \leq g(d) \leq \vse (c) \leq \vse( \nu ) \leq \vse(d) \leq h(c) \leq h(\nu) \leq h(d).
\]
By differentiation, it can be shown that for any fixed $\mu \in (0, 1)$, $\mscr{M}(z,\mu)$ is monotonically increasing with respect to $z \in
(0, \mu)$. Since $g (\nu) \leq g(d) \leq \vse(\nu)$ for all $\nu \in [c, d]$, it follows that \be \la{com1}
 \mscr{M} ( g (\nu), \vse(\nu) ) \leq
\mscr{M} ( g (d), \vse(\nu) ), \qqu \fa \nu \in [c, d]. \ee By differentiation, it can be shown that for any fixed $z \in (0, 1)$,
$\mscr{M}(z,\mu)$ is monotonically decreasing with respect to $\mu \in (z, 1)$. Since $g (d) \leq \vse(c) \leq \vse(\nu) \leq 1$ for all $\nu
\in [c, d]$, we have \be \la{com2}
 \mscr{M} ( g (d), \vse(\nu) ) \leq \mscr{M} ( g (d), \vse(c) ), \qqu \fa \nu \in [c, d].
\ee By virtue of (\ref{com1}) and  (\ref{com2}), we have \be \la{com3}
 \mscr{M} ( g (\nu), \vse(\nu) ) \leq \mscr{M} ( g (d), \vse(c) ), \qqu \fa \nu
\in [c, d]. \ee
Similarly, it can be shown that \bel &  & \mscr{M} ( h (\nu), \vse(\nu) ) \leq \mscr{M} ( h (c), \vse(d) ), \la{com4}\\
&  & \mscr{M} ( g (\nu), \vse(\nu) ) \geq \mscr{M} ( g (c), \vse(d) ), \la{com5}\\
&  & \mscr{M} ( h (\nu), \vse(\nu) ) \leq \mscr{M} ( h (d), \vse(c) ) \la{com6} \eel for all $\nu \in [c, d]$.  Combining (\ref{com3}),
(\ref{com4}), (\ref{com5}) and (\ref{com6}) yields (\ref{bounda}) and (\ref{boundb}).  Theorem \ref{Bounded_Mean_mix_general_Massart} is thus
established.

\end{document}